\documentclass[12pt]{amsart}
\usepackage{amssymb, amscd, amsmath, amsthm, latexsym, enumerate, eepic}

\newtheorem*{thm}{Theorem}

\begin{document}

\title{flat 2-orbifolds and seifert fibred 4-manifolds}

\author{J. A. Hillman }
\address{{School of Mathematics and Statistics, University of Sydney,}
\newline
Sydney,  NSW 2006, Australia }

\email{jonathan.hillman@sydney.edu.au}

\begin{abstract}
This is a summary of some of the basic facts about 
flat 2-orbifold groups, otherwise known as 
2-dimensional crystallographic groups.
We relate the geometric and topological presentations of these groups, 
and consider structures corresponding to decompositions of the orbifolds 
as fibrations or as unions.
We also consider  covering relations,
and record the bases of Seifert fibrations of 4-manifolds with geometries of 
solvable Lie type.
\end{abstract}



\maketitle

An $n$-dimensional crystallographic group is a discrete subgroup $\pi$ 
of the group $E(n)=\mathbb{R}^n\rtimes{O(n)}$ of isometries 
of euclidean $n$-space $\mathbb{R}^n$
which acts properly discontinuously on $\mathbb{R}^n$.
The {\it translation subgroup\/} $\pi\cap\mathbb{R}^n$ is  
a lattice of rank $n$, 
with quotient a finite subgroup of $O(n)$, 
called the {\it holonomy group\/} of $\pi$.
Since conjugation by $\pi$ preserves the lattice,
the holonomy group is conjugate in $GL(n,\mathbb{R})$ 
to a subgroup of $GL(n,\mathbb{Z})$.
The quotient $B$ of $\mathbb{R}^n$ by the action 
has a natural orbifold structure, 
recording the images of points with nontrivial stabilizers.
The group $\pi$ is then the {\it orbifold fundamental group\/}  
$\pi^{orb}(B)$.

It is well known that when $n=2$ there are just 17 possibilities.
We shall relate the presentations of the groups deriving from 
their structure as an extension of a finite group by a lattice 
to those deriving from the orbifold structure.
We give explicit embeddings of each group in $E(2)$, 
where there is not an obvious choice.
(However, we do not consider the issue of moduli, i.e., 
the parametrization of all such embeddings of a given flat 2-orbifold group.)
The orbifold fibres over a 1-orbifold if and only if 
the group is an extension of $\mathbb{Z}$ or  
the infinite dihedral group $D_\infty$.
In the latter case the group is also a generalized free product 
with amalgamation (GFPA), 
corresponding to a decomposition of the orbifold 
along a codimension-1 suborbifold.
In \S4 we describe the minimal proper covering relations 
between these orbifolds.
Finally we consider which flat 2-orbifolds can be bases of Seifert fibrations of 3- or 4-manifolds.
(This was the instigation for this work.)
Most of this material is well known; the main novelty here is 
perhaps in bringing this material together.
We include several relevant expository articles in the bibliography, 
although these are not explicitly invoked here.

\section{the holonomy group}

A finite subgroup of $SL(2,\mathbb{R})$ is conjugate into a subgroup of $SO(2)\cong{S^1}$,
and hence is cyclic.
Let 
\[
A=\left(\begin{smallmatrix}
0&1\\ -1&0
\end{smallmatrix}\right),\quad 
{B}=\left(\begin{smallmatrix}
0&-1\\ 1&1
\end{smallmatrix}\right)\quad\mathrm{and}\quad{R}=\left(\begin{smallmatrix}
0&1\\ 1&0
\end{smallmatrix}\right).
\] 
If $M$ in $GL(2,\mathbb{Z})$ has finite order $n$ then $\Delta_M=\det(tI-M)$
is quadratic, and its irreducible factors divide $t^n-1$.
Hence $n=1,2,3,4$ or 6.
If $n=3,4$ or 6 then $\Delta_M$ is irreducible and $\mathbb{Z}[t]/(\Delta_M)$
is the ring  $\mathbb{Z}[\zeta_n]$ of $n$th roots of unity, and is a PID.
It follows easily that $M$ is conjugate to $B^2$, $A$ or $B$, respectively.
Similarly, if $n=2$ then either $M=-I$ or $M$ is conjugate to $R$ or $AR$.
With a little calculation we see that if $G$ is a nontrivial finite subgroup of $GL(2,\mathbb{Z})$
it is conjugate to one of the cyclic groups generated by $A$, $-I$, $B$, $B^2$, $R$ or $AR$, 
or to a dihedral subgroup generated by $\{ A,R\}$, 
$\{-I, R\}$, $\{ -I,AR\}$, $\{ B,R\}$, $\{ B^2,R\}$ or $\{ B^2,RB\}$.

Let $\widetilde{\mathbb{Z}}^2$ denote $\mathbb{Z}^2$, considered as a $G$-module
via the inclusion $G<GL(2,\mathbb{Z})$.
The extensions of $G$ by $\widetilde{\mathbb{Z}}^2$ corresponding to  2-dimensional crystallographic group $\pi$ with holonomy $G$ are determined by  $H^2(G;\widetilde{\mathbb{Z}}^2)$.
In 10 of the 13 cases  $H^2(G;\widetilde{\mathbb{Z}}^2)=0$,
and so the semidirect product $\mathbb{Z}^2\rtimes{G}$
is the unique extension of $G$ by $\mathbb{Z}^2$ (up to isomorphism).
The eight semidirect products with 
$G\leq\langle{A,R}\rangle=GL(2,\mathbb{Z})\cap{O(2)}$
embed as discrete cocompact subgroups of $E(2)=\mathbb{R}^2\rtimes{O(2)}$
in the obvious way.
The five other semidirect products embed in 
$Aff(2)=\mathbb{R}^2\rtimes{GL(2,\mathbb{R})}$,
and these embeddings may be conjugated into $E(2)$.
There are four isomorphism classes of groups which are
non-split extensions, and we shall give explicit embeddings in each case.

Alternatively, we may compare local and global calculations of the orbifold Euler characteristic.
If $B$ is a flat 2-orbifold then $\chi^{orb}(B)=0$, since $B$ is finitely covered by the torus.
On the other hand, 
if $B$ is orientable then $B=F(a_1,\dots,a_k)$ for some closed surface $F$ 
and integers $a_i>1$,
and so 
\[
\chi^{orb}(B)=\chi(F)-k+\Sigma\frac1{a_i}.
\]
Hence either $F=T$ and $k=0$ or $F=S^2$ and $\Sigma\frac1{a_i}=k-2$.
A little arithmetic then shows that there are just five possibilities.
A similar calculation leads to the twelve non-orientable flat 2-orbifolds.
(See \cite{Sco}.
We may instead argue that if an orientation reversing  orbifold involution acts freely 
on the underlying surface $F$ then the quotient orbifold must be either $Kb$ or $P(2,2)$.
The remaining non-orientable flat 2-orbifolds are quotients by reflections across simple closed curves in $F$, fixing cone points on the reflector curve and permuting the others.)
We again find that there are 17 flat 2-orbifold groups.

\section{presentations}

Each group is identified by the traditional crystallographic symbol 
(in square brackets) and by the now-standard orbifold symbol
\cite{Mo}.
In each case we give first a presentation arising from the extension
and then one deriving from the corresponding flat orbifold.
Epimorphisms to $D_\infty$ correspond to GFPA structures, 
arising naturally from Van Kampen's Theorem.
These are essentially unique for $\mathbb{A}$, $\mathbb{M}b$ and $Kb$, 
since the kernel must contain the centre $\zeta\pi$.

Let $\mathbb{I}=[[0,1]]$ and $\mathbb{J}=[[0,1]$ be the reflector interval
and the interval with one reflector endpoint and one ordinary endpoint.
Let $Mb$ and $D(2,2)$ be the M\"obius band and the disc with two cone points, 
but with ordinary boundaries.
Then $\mathbb{I}$ is the one-point union of two copies of $\mathbb{J}$, 
and so $\pi^{orb}(\mathbb{J})=\mathbb{Z}/2\mathbb{Z}$ and
$\pi^{orb}(\mathbb{I})\cong\pi^{orb}(D(2,2))\cong{D_\infty}$. 

The first four orbifolds ($T,\mathbb{A},Kb$ and $\mathbb{M}b$) fibre over $S^1$.

\smallskip
{\bf Holonomy} $G=1$. 

\smallskip
\noindent$[p1]=T$. \quad $\mathbb{Z}^2=\langle{x,~y}\mid{xy=yx}\rangle$.  

In the subsequent presentations the generators $a,b,c,d,j,n$ and $r$ shall
represent elements whose images in the holonomy group have matrices $A,B,B^2,AR,-I, BR$ and $R$, respectively,
with respect to the basis $\{x,y\}$ for the translation subgroup $\mathbb{Z}^2$.
(The other generators $m,p,s,t,u,v,w,z$ do not have such fixed 
interpretations.)

\smallskip
{\bf Holonomy} $\mathbb{Z}/2\mathbb{Z}=\langle{AR}\rangle$.
In this case $H^2(G;\widetilde{\mathbb{Z}}^2)=\mathbb{Z}/2\mathbb{Z}$.

\smallskip
\noindent$[pm]=\mathbb{A}=S^1\times\mathbb{I}$. 
\quad $\pi=\mathbb{Z}\times{D_\infty}\cong
(\mathbb{Z}\oplus{\mathbb{Z}/2\mathbb{Z}})*_\mathbb{Z}(\mathbb{Z}\oplus{\mathbb{Z}/2\mathbb{Z}})$. 

This is the split extension.
\[\langle{\mathbb{Z}^2,~d}\mid{dx=xd,~dyd=y^{-1},~d^2=1}\rangle.\]
Let $u=dy$. Then also 
\[\langle{d,~u,~x}\mid{dx=xd,~ux=xu,~d^2=u^2=1}\rangle.\]

The subgroups $\langle{x+y,d}\rangle$ and $\langle{dx,y}\rangle$ 
are isomorphic to $\pi^{orb}(\mathbb{M}b)$ and $\pi_1(Kb)$, respectively.

\smallskip
\noindent$[pg]=Kb=Mb\cup{Mb}$. \quad $\pi=
\mathbb{Z}\rtimes_{-1}\mathbb{Z}\cong\mathbb{Z}*_\mathbb{Z}\mathbb{Z}$. 

This is the non-split extension.
\[\langle{\mathbb{Z}^2,~d}\mid{d^2=x,~dyd^{-1}=y^{-1}}\rangle=
\langle{d,~y}\mid{dyd^{-1}=y^{-1}}\rangle.\] 
Let $u=dy$. Then also 
\[\langle{d,~u}\mid{d^2=u^2}\rangle.\]
The quotient of $\pi$ by its centre $\langle{x}\rangle\cong\mathbb{Z}$ is $D_\infty$, 
but this extension does not split.

We may embed $\pi$ in $E(2)$ via 
${y}\mapsto(\mathbf{j},I_2)$ and
$z\mapsto(\frac12\mathbf{i},\left(\begin{smallmatrix}
1&0\\ 0&-1\end{smallmatrix}\right))$. 

\smallskip
{\bf Holonomy} $\mathbb{Z}/2\mathbb{Z}=\langle{R}\rangle$: 

\smallskip
\noindent$[cm]=\mathbb{M}b=Mb\cup{S^1\times\mathbb{J}}$.  \quad
$\pi\cong{\mathbb{Z}*_\mathbb{Z}(\mathbb{Z}\oplus\mathbb{Z}/2\mathbb{Z})}$.
\[\langle{\mathbb{Z}^2,~r}\mid{rxr=y,~r^2=1}\rangle.\]
Let $z=xr$. Then also 
\[\langle{r,~z}\mid{rz^2=z^2r,~r^2=1}\rangle.\]
Let $\tau$ be the involution of the normal subgroup 
$\langle{xy^{-1},r}\rangle=\langle{r,zrz^{-1}}\rangle\cong{D_\infty}$
which swaps $r$ and $zrz^{-1}$.
Thus $\pi\cong{D_\infty\rtimes}_\tau\mathbb{Z}$.
The centre is $\langle{xy}\rangle\cong\mathbb{Z}$ and $\pi/\zeta\pi\cong{D_\infty}$, 
but this extension does not split.

The subgroup $\langle{xr,rx}\rangle=\langle{z,[r,z]}\rangle$ 
is isomorphic to $\pi_1(Kb)$.

\medskip
The next five orbifolds fibre over $\mathbb{I}$.

\smallskip
{\bf Holonomy} $\mathbb{Z}/2\mathbb{Z}=\langle-I\rangle$: 

\smallskip
\noindent $[p2]=S(2,2,2,2)=D(2,2)\cup{D(2,2)}$.\quad 
$\pi\cong{D_\infty*_\mathbb{Z}D_\infty}$. 
\[\langle{\mathbb{Z}^2,~j}\mid{jxj=x^{-1},~jyj=y^{-1},~j^2=1}\rangle.\]
Let $u=jx$ and $v=jy$. 
Then also 
\[\langle{j,~u,~v}\mid{j^2=u^2=v^2=(juv)^2=1}\rangle.\]
This group is a semidirect product of the normal subgroup 
$\langle{ju}\rangle$ with $\langle{j,v}\rangle$.
Thus $\pi\cong\mathbb{Z}\rtimes{D_\infty}$.

\smallskip
{\bf Holonomy} $D_4=\langle-I,R\rangle$: 

\smallskip
\noindent$[cmm]=\mathbb{D}(2,\overline{2},\overline{2})$. \quad 
$\pi\cong\langle{S(2,2,2,2),~r}\mid{rxr=y,~r^2=(jr)^2=1}\rangle$.
Let $u=jx$, $v=jy$ and $z=jr$. 
Then $v=rur$ and $juv=rzurur$. Then also
\[\langle{r,~u,~z}\mid{r^2=u^2=z^2=(rz)^2=(zuru)^2=1}\rangle.\]
This group is the semidirect product of the normal subgroup
$\langle{xy^{-1}\!,u}\rangle$ with $\langle{j,x}\rangle$,
and so $\pi\cong{D_\infty\rtimes{D_\infty}}$.
The corresponding GFPA structure
$\pi\cong(D_\infty\times{\mathbb{Z}/2\mathbb{Z}})*_{D_\infty}D_\infty$ 
derives from the decomposition of the disc along a chord which
separates the cone point from the corner points.

\smallskip
{\bf Holonomy}  $D_4=\langle-I,AR\rangle$. 
In this case $H^2(G;\widetilde{\mathbb{Z}}^2)=(\mathbb{Z}/2\mathbb{Z})^2$.

\noindent{The group $\pi$ has a presentation}
\[
\langle{\mathbb{Z}^2,~d,~j}\mid
{dx=xd,~dyd^{-1}=y^{-1},~jxj=x^{-1},~jyj=y^{-1},}
\]
\[(jd)^2=y^e,~d^2=x^f,~j^2=1\rangle.\]
We may assume that $0\leq{e,f}\leq1$.
In all cases, $\langle{x}\rangle$ and $\langle{y}\rangle$ 
are the maximal infinite cyclic normal subgroups.

Two extension classes give isomorphic groups, and so
there are three possibilities:

\smallskip
\noindent
$[pmm]=\mathbb{D}(\overline{2},\overline{2},\overline{2},\overline{2})$.\quad
$\pi\cong(D_\infty\times{\mathbb{Z}/2\mathbb{Z}})*_{D_\infty}(D_\infty\times\mathbb{Z}/2\mathbb{Z})$.

This is the split extension (with $e=f=0$).
It is also $D_\infty\times{D_\infty}$:
\[\langle{\mathbb{Z}^2,~d,~j}
\mid{dx=xd,~dyd=y^{-1},~jxj=x^{-1},~jyj=y^{-1},}\] 
\[{d^2=j^2=(dj)^2=1}\rangle.\]
Let $s=jdx$ and $t=dy$. Then also 
\[\langle{d,~j,~s,~t}\mid
{d^2=j^2=s^2=t^2=(st)^2=(tj)^2=(jd)^2=(ds)^2=1}\rangle,\]
\[
\mathrm{or}\quad
\langle{d,~j,~s,~t}\mid{d,t\leftrightharpoons{j,s},~d^2=j^2=s^2=t^2=1}\rangle.\]
The GFPA structure derives from the decomposition of the disc 
along a chord which separates one pair of adjacent corner points 
from the others.

\smallskip
\noindent$[pmg]=\mathbb{D}(2,2)=S^1\times\mathbb{J}\cup{D(2,2)}$.\quad 
$\pi\cong(\mathbb{Z}\oplus\mathbb{Z}/2\mathbb{Z})*_\mathbb{Z}D_\infty$.

This corresponds to $(e,f)=(1,0)$.
(The choice $(0,1)$ gives an isomorphic group.)
\[\langle{\mathbb{Z}^2,~d,~j}\mid{jxj=x^{-1},
~y=(jd)^2,~dx=xd,~d^2=j^2=1}\rangle.\]
Let $v=jx$. Then also
\[\langle{d,~j,~v}\mid{djv=jvd,~d^2=j^2=v^2=1}\rangle.\]
The relation $djv=jvd$ is equivalent to $jdj=vdv$, since 
${j^2=v^2=1}$.
Hence we also have $\pi\cong{D_\infty*_{D_\infty}D_\infty}$, 
This GFPA structure derives from the decomposition
of the disc along a chord which separates the two cone points.

The subgroups $\langle{jv}\rangle$ and $\langle{d,jdj}\rangle$ are normal, 
and so $\pi\cong\mathbb{Z}\rtimes{D_\infty}$ and 
$\pi\cong{D_\infty\rtimes{D_\infty}}$. 

We may embed $\pi$ in $E(2)$ via 
$d\mapsto(-\frac12\mathbf{j},\left(\begin{smallmatrix}
1&0\\ 0&-1\end{smallmatrix}\right))$, $j\mapsto(0,-I_2)$
and ${v}\mapsto(-\mathbf{i},-I_2).$

The index-2 subgroups of the group of $\mathbb{D}(2,2)$ 
corresponding to $\langle{AR}\rangle<D_4$ and $\langle{-AR}\rangle<D_4$
are isomorphic to 
$\pi_1(Kb)=\mathbb{Z}\rtimes_{-1}\mathbb{Z}$ and $\pi^{orb}(\mathbb{A})=\mathbb{Z}\times{D_\infty}$, respectively.

\smallskip
\noindent$[pgg]=P(2,2)=Mb\cup{D(2,2)}$.\quad
$\pi\cong\mathbb{Z}*_\mathbb{Z}D_\infty$.
 
This corresponds to $e=f=1$.
\[\langle{\mathbb{Z}^2,~d,~j}\mid{d^2=x,~(jd)^2=y,~jd^2j=d^{-2},~
d^2(jd)^2=(jd)^2d^2,~j^2=1}\rangle.\]
Let $v=jd^2$. Then this reduces to
\[
\langle{d,~j,~v}\mid{d^2=jv,~j^2=v^2=1}\rangle\quad\mathrm{or~just}\quad
\langle{d,~j}\mid{(jd^2)^2=j^2=1}\rangle.
\]
There is an automorphism which fixes $j$ and swaps $d$ and $jd$.
We have $\pi/\langle{(jd)^2}\rangle\cong\pi/\langle{d^2}\rangle\cong{D_\infty}$,
but these extensions do not split.

We may embed $\pi$ in $E(2)$ via 
$d\mapsto(\frac12\mathbf{i},\left(\begin{smallmatrix}
1&0\\ 0&-1\end{smallmatrix}\right))$ and
$j\mapsto(\frac12(\mathbf{i}+\mathbf{j}),-I_2)$.

The index-2 subgroups of the group of $P(2,2)$
corresponding to $\langle{AR}\rangle<D_4$ and $\langle{-AR}\rangle<D_4$
are both isomorphic to $\pi_1(Kb)$.
Hence non-trivial elements of finite order in $\pi$ must have order 2, 
and be of the form $wj$, with $w\in\langle{d,x,y}\rangle$.
In fact we must have $w\in\mathbb{Z}^2$, and so
all such elements are in
the orientation-preserving subgroup 
$\pi^+=\langle{j,x,y}\rangle=\pi^{orb}(S(2,2,2,2))$.

The remaining eight orbifolds do not fibre over $S^1$ or $\mathbb{I}$.

\smallskip
{\bf Holonomy} $\mathbb{Z}/4\mathbb{Z}=\langle{A}\rangle$: 

\smallskip
\noindent$[p4]=S(2,4,4)$. \quad
$\langle{\mathbb{Z}^2,~a}\mid{axa^{-1}=y^{-1},~aya^{-1}=x,~a^4=1}\rangle$.

Let $j=a^2x$. Then also $\langle{a,j}\mid{a^4=j^2=(aj)^4=1}\rangle$.

\smallskip
{\bf Holonomy} $D_8=\langle{A,R}\rangle$.\quad 
In this case $H^2(G;\widetilde{\mathbb{Z}}^2)=\mathbb{Z}/2\mathbb{Z}$.

\noindent{The group $\pi$ has a presentation}
\[\langle{\mathbb{Z}^2,~a,~r}\mid{axa^{-1}=y^{-1},~aya^{-1}=x,~rxr=y,~a^4=r^2=1,~
(ar)^2=y^e}\rangle,\]
where $e=0$ or 1.
Let $j=a^2x$, so $x=a^2j$ and $y=aja$.
Then this reduces to
\[
\langle{a,j,r}\mid{a^4=j^2=(aj)^4=r^2=1,~ra^2jr=aja,~(ar)^2=(aja)^e}\rangle.
\]

\noindent$[p4m]=\mathbb{D}(\overline{2},\overline{4},\overline{4})$.\quad

This is the split extension (with $e=0$).
It reduces to
\[\langle{a,j,r}\mid{a^4=j^2=(aj)^4=r^2=(ar)^2=1,~arj=jar}\rangle,\]
since $a^2r=ra^2$ when $e=0$.
Let $w=ar$ and $z=jw=jar$. Then also
\[\langle{r,w,z}\mid{r^2=w^2=z^2=(rz)^4=(zw)^2=(wr)^4=1}\rangle.\] 

\noindent$[p4g]=\mathbb{D}(\overline{2},4)$.

This is the non-split extension (with $e=1$).
\[\langle{a,j,r}\mid{a^4=j^2=(aj)^4=r^2=1,~ra^2jr=aja,~ (ar)^2=aja}\rangle.
\]
Hence $aj=a^{-1}rara=(ra)^{-1}ara$ and $j=rara^{-1}$, 
so this reduces to 
\[\langle{a,r}\mid{a^4=r^2=(rara^{-1})^2=1}\rangle.\]

We may embed $\pi$ in $E(2)$ via 
$a\mapsto(0,\left(\begin{smallmatrix}
0&-1\\1&0\end{smallmatrix}\right))$,
$j\mapsto(-\mathbf{i},-I_2)$ and
${r}\mapsto(\frac12\mathbf{i},\left(\begin{smallmatrix}
0&1\\ 1&0\end{smallmatrix}\right))$.

The index-2 subgroup of this group corresponding to 
$\langle-I,AR\rangle<D_8$ is the split extension.
In particular, the groups with holonomy $\mathbb{Z}/4\mathbb{Z}$ or $D_8$ do not contain
the groups of $Kb$, $\mathbb{D}(2,2)$ or $P(2,2)$.

\smallskip
{\bf Holonomy} $\mathbb{Z}/3\mathbb{Z}=\langle{B^2}\rangle$: 

\smallskip
\noindent$[p3]=S(3,3,3)$. \quad
$\langle{\mathbb{Z}^2,c}\mid{cxc^{-1}=x^{-1}y,~cyc^{-1}=x^{-1},~c^3=1}\rangle$.

Let $u=cx$. Then also 
\[\langle{c,u}\mid{c^3=u^3=(cu)^3=1}\rangle.\]

{\bf Holonomy} $\mathbb{Z}/6\mathbb{Z}=\langle{B}\rangle$: 

\smallskip
\noindent $[p6]=S(2,3,6)$.\quad
$\langle{\mathbb{Z}^2,~b}\mid{bxb^{-1}=y,~byb^{-1}=x^{-1}y,~b^6=1}\rangle$.

Let $v=b^2x$. Then also \[\langle{b,v}\mid{b^6=v^3=(bv)^2=1}\rangle.\] 

{\bf Holonomy} $D_6=\langle{B^2,R}\rangle$: 

\smallskip
\noindent$[p3m1]=\mathbb{D}(\overline{3},\overline{3},\overline{3})$.\quad
$\langle{S(3,3,3),~r}\mid{rxr=y,~r^2=(rc)^2=1}\rangle$.

Let  $s=rc$ and $t=crx$. Then also 
\[\langle{r,s,t}\mid{r^2=s^2=t^2=(rs)^3=(st)^3=(tr)^3=1}\rangle.\]

{\bf Holonomy} $D_6=\langle{B^2,BR}\rangle$: 

\smallskip
\noindent$[p31m]=\mathbb{D}(3,\overline{3})$. \quad
$\langle{S(3,3,3),~n}\mid{nxn=x^{-1}y,~ny=yn,~n^2=(nc)^2=1}\rangle$.

Let $v=nx^{-1}y^{-1}$ and $w=nx^2y^{-1}$. Then also
\[
\langle{c,v,w}\mid{c^3=v^2=w^2=1,~w=cvc^{-1}}\rangle
\]
\[
\mathrm{or}\quad
\langle{c,w}\mid{w^2=c^3=(c^{-1}wcw)^3=1}\rangle.
\]

{\bf Holonomy} $D_{12}=\langle{B,R}\rangle$: 

\smallskip
\noindent$[p6m]=\mathbb{D}(\overline{2},\overline{3},\overline{6})$.\quad
$\langle{S(2,3,6),~r}\mid{rxr=y,~r^2=(rb)^2=1}\rangle$.

Let $m=brb^{-1}=b^2r$, $n=br$ and $p=rbxy^{-2}$. Then also
\[\langle{m,n,p}\mid{m^2=n^2=p^2=(mn)^6=(np)^3=(pm)^2=1}\rangle.\] 

The matrix $B$ is not orthogonal.
However conjugation by $\left(\begin{smallmatrix}
-2&1\\ 0&\sqrt3\end{smallmatrix}\right))$
carries $\langle{B,R}\rangle$ into $O(2)$,
and thus carries each of the subgroups $\mathbb{Z}^2\rtimes{G}$ of $Aff(2)$ 
determined by the five groups with 3-torsion into $E(2)$.

In practice it can be useful to first sort these groups by their abelianizations,
which are:

$\mathbb{Z}^2$ for $T$, \quad 

$\mathbb{Z}\oplus(\mathbb{Z}/2\mathbb{Z})^2$ for $\mathbb{A}$,

$\mathbb{Z}\oplus{\mathbb{Z}/2\mathbb{Z}}$ for $Kb$ and $\mathbb{M}b$,\quad

$(\mathbb{Z}/2\mathbb{Z})^4$ for $\mathbb{D}(\overline2,\overline2,\overline2,\overline2)$,

$(\mathbb{Z}/2\mathbb{Z})^3$ for $S(2,2,2,2)$, $\mathbb{D}(2,2)$,
$\mathbb{D}(2,\overline2,\overline2)$ and $\mathbb{D}(\overline2,\overline4,\overline4)$,

$(\mathbb{Z}/2\mathbb{Z})^2$ for $\mathbb{D}(\overline2,\overline3,\overline6)$,\quad

$\mathbb{Z}/2\mathbb{Z}$ for $\mathbb{D}(\overline3,\overline3,\overline3)$,

$\mathbb{Z}/4\mathbb{Z}\oplus{\mathbb{Z}/2\mathbb{Z}}$ for $S(2,4,4)$, $P(2,2)$ and 
$\mathbb{D}(\overline2,4)$,

$\mathbb{Z}/3\mathbb{Z}$ for $S(3,3,3)$\quad{and}\quad  

$\mathbb{Z}/6\mathbb{Z}$ for $S(2,3,6)$ and $\mathbb{D}(3,\overline3)$.

Only the groups of the pairs $\{Kb,\mathbb{M}b\}$ and $\{\mathbb{D}(2,2),\mathbb{D}(2,\overline2,\overline2)\}$
have both isomorphic holonomy groups and isomorphic abelianizations, 
and in each case one of the pair is a semidirect product.

All but one of the groups with finite abelianization are generated by elements of finite order,
the exception being $\pi^{orb}(P(2,2))$.

\section{fibrations over $\mathbb{I}$}

We shall say that two epimorphisms 
$\lambda,\lambda':\pi\to{D_\infty}=\pi^{orb}(\mathbb{I})$ 
are equivalent if $\lambda'=\delta\lambda\alpha$ for some $\alpha\in{Aut(\pi)}$
and $\delta\in{Aut(D_\infty)}$.
Such epis are most easily found by considering the possible kernels,
which are maximal normal virtually-$\mathbb{Z}$ subgroups.
(Note that since $Out(D_\infty)=\mathbb{Z}/2\mathbb{Z}$ and we are working with epimorphisms, 
it is sufficient to take $\delta$ to be either the automorphism which swaps 
the standard generators or the identity.) 

All epimorphisms from $\pi^{orb}(S(2,2,2,2))$ to $D_\infty$ are equivalent.
These correspond to fibrations over $\mathbb{I}$ with general fibre $S^1$
and two singular fibres (reflector intervals connecting pairs of cone points).

All epimorphisms from 
$\pi^{orb}(\mathbb{D}(\overline{2},\overline{2},\overline{2},\overline{2}))$
to $D_\infty$ are equivalent.
(There are two maximal normal virtually-$\mathbb{Z}$ subgroups.)
These correspond to the projections of
$\mathbb{D}(\overline{2},\overline{2},\overline{2},\overline{2})=
\mathbb{I}\times\mathbb{I}$ onto its factors.

All epimorphisms from $\pi^{orb}(P(2,2))$ to $D_\infty$ are equivalent.
(There are two maximal normal virtually-$\mathbb{Z}$ subgroups.)
These correspond to a fibration over $\mathbb{I}$ with general fibre $S^1$
and two singular fibres (the centreline of $Mb$,
and one reflector interval connecting the cone points).

There are two equivalence classes of epimorphisms from 
$\pi^{orb}(\mathbb{D}(2,2))$ to $D_\infty$.
(There are two maximal normal virtually-$\mathbb{Z}$ subgroups.)
One corresponds to a fibration over $\mathbb{I}$ with general fibre $S^1$
and two singular fibres (the reflector curve and a reflector interval
connecting the cone points).
The other corresponds to a fibration with general fibre $\mathbb{I}$ 
and two singular fibres (two reflector intervals,
each connecting a cone point to a reflector curve).

\setlength{\unitlength}{1mm}
\begin{picture}(90,42)(-75,-2)

\put(-45,22){\circle{28}}
\put(-45,22){\circle{26}}

\put(-53,21){$\bullet$}
\put(-39,21){$\bullet$}
\qbezier(-52,22)(-45,22)(-38,22)

\qbezier(-50,19)(-62,22)(-50,25)
\qbezier(-40,19)(-28,22)(-40,25)
\qbezier(-50,25)(-45,26)(-40,25)
\qbezier(-50,19)(-45,18)(-40,19)

\put(-51,23){$2$}
\put(-37,23){$2$}

\put(10,22){\circle{28}}
\put(10,22){\circle{26}}
\put(2,21){$\bullet$}
\put(16,21){$\bullet$}

\put(4,23){$2$}
\put(18,23){$2$}
\qbezier(-4,22)(-1,22)(3,22)
\qbezier(17,22)(21,22)(24,22)
\qbezier(10,8)(10,22)(10,36)
\qbezier(5,9)(10,22)(5,35)
\qbezier(-1.5,14)(10,22)(-1.5,30)

\qbezier(15,9)(10,22)(15,35)
\qbezier(21.5,14)(10,22)(21.5,30)

\put(-25,0){$\mathbb{D}(2,2)$}

\end{picture}

All epimorphisms from $\pi^{orb}(\mathbb{D}(2,\overline{2},\overline{2}))$
to $D_\infty$ are equivalent.
(There are two maximal normal virtually-$\mathbb{Z}$ subgroups.)
These correspond to a fibration over $\mathbb{I}$ 
with general fibre $\mathbb{I}$ and one exceptional fibre
(a reflector interval connecting the cone point to a reflector curve).

\setlength{\unitlength}{1mm}
\begin{picture}(90,43)(2,-1)

\qbezier(35,22)(60,50)(85,22)
\qbezier(35,22)(60,-6)(85,22)
\qbezier(35.5,21.5)(60,49)(84.5,21.5)
\qbezier(35.5,22.5)(60,-5)(84.5,22.5)
\put(59,23){$2$}

\qbezier(38,19)(60,48)(82,19)

\qbezier(45,13)(60,52)(75,13)
\qbezier(50,10.5)(60,50)(70,10.5)
\qbezier(55,8.5)(60,48)(65,8.5)
\qbezier(60,22)(60,15)(60,8)
\put(59,21){$\bullet$}
\put(52,0){$\mathbb{D}(2,\overline{2},\overline{2})$}

\end{picture}

\section{coverings}

If $\alpha$ and $\beta$ are two flat orbifold groups 
and there is a monomorphism $\alpha\to\beta$ which
is an isomorphism on the translation subgroups then
we shall say that the corresponding orbifold cover is {\it equitranslational}.
In this case $\beta_1(\alpha)\geq\beta_1(\beta)$,
the holonomy group $G_\alpha$ is conjugate in $GL(2,\mathbb{Z})$ 
to a subgroup of $G_\beta$,
and the extension class $e_\beta\in{H^2(G_\beta;\widetilde{\mathbb{Z}}^2)}$ 
must restrict to $e_\alpha$.

Such inclusions are easily determined from the lattices 
of subgroups of the two maximal subgroups $\langle{A,R}\rangle$ and
$\langle{B,R}\rangle$ of $GL(2,\mathbb{Z})$ (modulo conjugacy).
Determining the lattices of equitranslational covers of
$\mathbb{D}(\overline2,\overline4,\overline4)$ or of
$\mathbb{D}(\overline2,\overline3,\overline6)$ is straightforward,
since in these cases the orbifold groups arising are split extensions.
The only subtle point is how the nonsplit extensions (with holonomy a 2-group)
restrict over subgroups of the holonomy groups.
The index-2 subgroups of $\pi^{orb}(\mathbb{D}(\overline2,4))$ corresponding 
to $\langle -I,AR\rangle$ are split extensions, and so
$\mathbb{D}(2,2)$ and $P(2,2)$ do not cover $\mathbb{D}(\overline2,4)$.
The index-2 subgroups of $\pi^{orb}(P(2,2))$ corresponding 
to $\langle \pm{AR}\rangle$ are both $\pi_1(Kb)$, and so
$\mathbb{A}$ does not cover $P(2,2)$.

The orientable orbifold $S(2,2,2,2)$ covers every flat orbifold with 
holonomy of even order except for $Kb$, $\mathbb{A}$, $\mathbb{M}b$,
$\mathbb{D}(3,\overline{3})$ and
$\mathbb{D}(\overline{3},\overline{3},\overline{3})$.
In particular,
$\mathbb{D}(2,2)$, $\mathbb{D}(2,\overline{2},\overline{2})$ and
$\mathbb{D}(\overline{2},\overline{2},\overline{2},\overline{2})$
are quotients by reflections across circles 
passing through 0, 2 or 4 cone points, respectively. 
The quotient by rotation of order 2 with both fixed points cone points 
is $S(2,4,4)$, 
while the quotient by rotation of order 3 with one fixed point 
a cone point is $S(2,3,6)$.

Similarly, $S(2,4,4)$ covers each of $\mathbb{D}(\overline{2},4)$ 
and $\mathbb{D}(\overline{2},\overline{4},\overline{4})$,
and $S(3,3,3)$ covers every flat orbifold with holonomy divisible by 3.

If we drop the requirement that the translation subgroups coincide,
there also less obvious inclusions.
It remains necessary that $G_\alpha$ be conjugate in $GL(2,\mathbb{Q})$ 
to a subgroup of $G_\beta$.
We may also use the (non)existence of reflector curves
and/or corner points in testing whether one orbifold covers another.
In all cases the subgroup generated by $2x$ and $2y$ is normal 
and of index 4 in the translation subgroup, and so
the orbifolds have degree-4 self-coverings.
(If $G\leq\langle{A,R}\rangle$ there is a degree-2 self-covering, 
since the subgroup generated by $x+y$ and $x-y$ is normal 
and of index 2 in the translation subgroup.)

For example,
although $AR$ and $R$ are not conjugate in $GL(2,\mathbb{Z})$,
they are conjugate in $GL(2,\mathbb{Z}[\frac12])$.
The inclusion $\mathbb{Z}\times{D_\infty}<{D_\infty\rtimes_\tau\mathbb{Z}}$
corresponds to the geometric fact that $\mathbb{A}$ covers $\mathbb{M}b$.
Folding $\mathbb{M}b$ across its centerline gives $\mathbb{A}$, 
while the quotient of $Kb$ by fibrewise reflection is $\mathbb{M}b$.
However $P(2,2)$ is covered only by $T,Kb$, $S(2,2,2,2)$ and itself.

The involution $[x:y:z]\mapsto[x:-y:-z]$ of $RP^2$
has one fixed point and one fixed circle.
Hence $P(2,2)$ covers $\mathbb{D}(2,2)$.

Rotating $\mathbb{D}(\overline{2},\overline{2},\overline{2},\overline{2})$
about its centre gives $\mathbb{D}(2,\overline{2},\overline{2})$.
Conversely, folding $\mathbb{D}(2,\overline{2},\overline{2})$
across a diameter through the cone point and separating the corner points
gives $\mathbb{D}(\overline{2},\overline{2},\overline{2},\overline{2})$.

Folding $\mathbb{D}(2,2)$ across a diameter separating the cone points
gives $\mathbb{D}(2,\overline{2},\overline{2})$.
Since $\mathbb{D}(2,2)$ has no corner points, 
it is not covered by $\mathbb{D}(2,\overline{2},\overline{2})$.

Folding $\mathbb{D}(\overline{2},4)$ across 
a diameter through the cone and corner points gives 
$\mathbb{D}(\overline{2},\overline{4},\overline{4})$.
However, $\mathbb{D}(\overline{2},4)$ is not covered by 
$\mathbb{D}(\overline{2},\overline{4},\overline{4})$,
since it has no corner points with stabilizer of order $>4$.
(On the algebraic side,
$\pi^{orb}(\mathbb{D}(\overline{2},\overline{4},\overline{4}))\cong\mathbb{Z}^2\rtimes{D_8}$
is the split extension, whereas $\pi^{orb}(\mathbb{D}(\overline{2},4))$ is the non-split extension,
which does not have $D_8$ as a subgroup.)

Every flat orbifold with holonomy divisible by 3 covers
$\mathbb{D}(\overline{2},\overline{3},\overline{6})$.

Both $\pi^{orb}(\mathbb{D}(\overline{3},\overline{3},\overline{3}))$
and $\pi^{orb}(\mathbb{D}(3,\overline{3}))$ are split extensions, 
and each contains the other as a subgroup of index 3.
(The normal subgroup of index 3 in 
$\pi^{orb}(\mathbb{D}(3,\overline{3}))\cong\mathbb{Z}^2\rtimes{D_6}$ 
is isomorphic to 
$\pi^{orb}(\mathbb{D}(\overline{3},\overline{3},\overline{3}))$, 
and this in turn contains
the subgroup generated by the holonomy group $D_6$ and $3\mathbb{Z}^2$,
which is isomorphic to  $\pi^{orb}(\mathbb{D}(3,\overline{3}))$.
Rotating $\mathbb{D}(\overline{3},\overline{3},\overline{3})$
about its centre gives $\mathbb{D}(3,\overline{3})$.
Is there a simple geometric description of the irregular 3-fold cover of
$\mathbb{D}(\overline{3},\overline{3},\overline{3})$ by $\mathbb{D}(3,\overline{3})$?

\section{seifert fibred 3-manifolds}

A Seifert fibred 3-manifold $M$ is flat or is a $\mathbb{N}il^3$-manifold if and only if
it is Seifert fibred over a flat 2-orbifold $B$. 
The fibration derives from an $S^1$-action
(i.e., the image of the general fibre in $\pi_1(M)$ is central) 
if and only if the orientation character $w=w_1(M)$
factors through $\pi^{orb}(B)$.
Every flat 3-manifold and every $\mathbb{N}il^3$-manifold is Seifert fibred,
but some flat 3-manifolds have several such fibrations, with distinct bases.
However, not all flat 2-orbifolds arise as bases of such fibrations.

There are just ten flat 3-manifolds, and it is easy to check the possibilities by hand.
Let $M$ be a flat 3-manifold which is Seifert fibred over $B$,
and let $h\in\pi=\pi_1(M)$ generate the image of the fundamental group of the general fibre.
Then $\pi^{orb}(B)\cong\pi/\langle{h}\rangle$.
The image of the translation subgroup of $\pi$ in $\pi^{orb}(B)$ is an abelian normal subgroup.
It follows that the holonomy of $B$ is a quotient of the holonomy of $M$,
and thus is cyclic or $(\mathbb{Z}/2\mathbb{Z})^2$ (by the classification of flat 3-manifold groups).
On considering the possible infinite cyclic normal subgroups we find that the possible bases
are the seven flat 2-orbifolds with no reflector curves, together with $\mathbb{A}$,
$\mathbb{M}b$ and $\mathbb{D}(2,2)$.
(We may exclude $\mathbb{D}(\overline2, \overline2, \overline2, \overline2)$,
since every flat 3-manifold group can be generated by at most 3 elements, whereas
 $\beta_1(\pi^{orb}( \mathbb{D}(\overline2, \overline2, \overline2, \overline2));\mathbb{F}_2)=4$.
 Similarly, since $\mathbb{D}(2,\overline2, \overline2)$ has holonomy $(Z/2Z)^2$ and
 $\beta_1(\pi^{orb}( \mathbb{D}(2,\overline2, \overline2));\mathbb{F}_2)=3$,
 it could only be the base of a Seifert fibration of the non-orientable flat orbifold $B_3$.
 But this can be ruled out by inspection.)
 
In a little more detail: the 3-torus $G_1$ has base $T$, the half-turn 3-manifold $G_2$ 
has bases $Kb$ and $S(2,2,2,2)$, $G_3$ has base $S(3,3,3)$, $G_4$ has base $S(2,4,4)$,
$G_5$ has base $S(2,3,6)$ and $G_6$ has base $P(2,2)$.
The fibrations of $G_3,\dots,G_6$ are unique.
The non-orientable flat 3-manifolds all have Seifert fibrations with base $Kb$.
In addition, $B_1=Kb\times{S^1}$ also has bases $T$ and $\mathbb{A}$,
$B_2$ also has bases $T$ and $\mathbb{M}b$,
$B_3$ also has bases $\mathbb{A}$ and $\mathbb{D}(2,2)$,
and $B_4$ also has base $\mathbb{M}b$.

In the $\mathbb{N}il^3$ case,  the Seifert fibration is unique.
The images of elements of finite order in $\pi^{orb}(B)$
in the holonomy group have determinant $+1$, 
by part (3) of Theorem 8.6 of \cite{Hi02}.
Hence $B$ can have no reflector curves, and so must be one of
$T$, $Kb$, $S(2,2,2,2)$, $P(2,2)$, $S(2,4,4)$, $S(3,3,3)$ or $S(2,3,6)$.
Each of these is the base of a Seifert fibration of a $\mathbb{N}il^3$-manifold.
(See \cite{De}.)
Since $\mathbb{N}il^3$-manifolds are orientable,
the fibration derives from an $S^1$-action 
if and only if the base is orientable,
and so is one of $T$, $S(2,2,2,2)$, $S(2,4,4)$, $S(3,3,3)$ or $S(2,3,6)$.

\section{seifert fibred 4-manifolds}

In dimension 4, manifolds which are Seifert fibred over a flat 2-orbifold 
with general fibre a flat 2-manifold are also geometric.
(However there are non-geometric 4-manifolds which are Seifert fibred 
over hyperbolic base 2-orbifolds.)
The relevant geometries are $\mathbb{E}^4$,
$\mathbb{N}il^3\times\mathbb{E}^1$, 
$\mathbb{N}il^4$ and $\mathbb{S}ol^3\times\mathbb{E}^1$.

All but three of the 74 flat 4-manifolds have Seifert fibrations. 
(See Chapter 8 of \cite{Hi02}.)
If the holonomy $G$ of $\pi^{orb}(B)$ is cyclic
then $B$ is the base of a Seifert fibration of a flat 3-manifold $N$,
and so $N\times{S^1}$ is a flat 4-manifold which is Seifert fibred with general fibre 
$T$ and base $B$.
These products are also total spaces of $T$-bundles over $T$.
If the holonomy $G$ is dihedral it is a quotient of $\pi_1(Kb)$,
and pulling back the extension of $G$ by $\mathbb{Z}^2$ over $\pi_1(Kb)$
gives the group of a flat 4-manifold which is Seifert fibred over $B$.
These 4-manifolds are also total spaces of $T$-bundles over $Kb$.
A fibration with general fibre $T$ derives from a $T$-action
(i.e., the image of $\pi_1(T)$ in $\pi=\pi_1(M)$ is central) 
if and only if $w=w_1(\pi)$ factors through $\pi^{orb}(B)$.

Manifolds with one of the geometries 
$\mathbb{N}il^3\times\mathbb{E}^1$, 
$\mathbb{N}il^4$ or $\mathbb{S}ol^3\times\mathbb{E}^1$
have universal cover diffeomorphic to a connected 4-dimensional solvable Lie group.
These Lie groups have characteristic subgroups $\mathbb{R}^2$;
these are the centre $\zeta{Nil^3}\times\mathbb{R}$, the hypercentre $\zeta_2Nil^4$
and the commutator subgroup ${Sol^3}'$,
respectively.
The foliations by cosets of such subgroups induce canonical
Seifert fibrations on the quotient manifolds.
If the geometry is $\mathbb{N}il^3\times\mathbb{E}^1$ or
$\mathbb{N}il^4$ the manifold is an infranilmanifold,
and these are treated at length in \cite{De}.

If the geometry is $\mathbb{N}il^4$ or $\mathbb{S}ol^3\times\mathbb{E}^1$
the manifold has an unique Seifert fibration.
The canonical Seifert fibration of  a $\mathbb{N}il^3\times\mathbb{E}^1$-manifold
lifts to the orbit fibration of a $T$-action on a finite cover
(i.e., $\pi_1(F)$ is commensurable with the centre of $\pi$), 
but in general such manifolds have other Seifert fibrations.
(For instance, if $N$ is a $\mathbb{N}il^3$-manifold with $\beta_1(N)=2$
then $N\times{S^1}$ is the total space of a $T$-bundle over $T$ in at least two
essentially distinct ways.)

The general fibre $F$ of a Seifert fibration of a  $\mathbb{N}il^3\times\mathbb{E}^1$-manifold
can be either $T$ or $Kb$, 
and inspection of the tables in Chapter 7 of \cite {De} shows that every flat 2-orbifold is the basis of 
some such fibration with $F=T$.
However, it is not clear from a cursory inspection whether this is always so 
for the canonical Seifert fibration.

\begin{thm}
Let $M$ be a $\mathbb{N}il^3\times\mathbb{E}^1$-manifold which is Seifert fibred 
over $B$.
If the general fibre $F=T$ and $\pi_1(F)$ is central in $\pi=\pi_1(M)$ then $B$ is orientable.
If $F=Kb$ then $B$ is $T$, $Kb$, $S(2,2,2,2)$, $S(3,3,3)$, $S(2,3,6)$ or $P(2,2)$.
\end{thm}

\begin{proof}
Suppose first that $F=T$.
If $a,b\in\pi$ represent a basis of the translation subgroup of $\pi^{orb}(B)$
and $t\in\pi$ then $[a,b]$ is a nontrivial element of $\pi_1(F)$, and $t[a,b]t^{-1}=[a,b]^{\det(t)}$,
where $\det(t)$ is the determinant of the image of $t$ in the holonomy of $\pi^{orb}(B)$.
Therefore if the image of $\pi_1(F)$ in $\pi$ is central then $B$ must be orientable.

Suppose now that $F=Kb$, and let $\kappa=\pi_1(F)$ and $w=w_1(\pi)$.
Then $w|_\kappa=w_1(\kappa)$.
Hence the orientable double cover $M^+$ is Seifert fibred over $B$ with general fibre $F^+=T$.
If the image of $g\in\pi$ in $\beta=\pi^{orb}(B)$ has finite order $o(g)$ then the subgroup
$\langle\kappa,g\rangle$ is a finite torsion-free extension of $\kappa$, 
and so is again isomorphic to $\pi_1(Kb)$.
If $x,y$ are generators for $\kappa$ such that $xyx^{-1}=y^{-1}$ then
$gxg^{-1}=x^\epsilon{y^s}$, for some $\epsilon=\pm1$ and $s\in\mathbb{Z}$,
and $gyg^{-1}=y^\eta$, for some $\eta=\pm1$.
Since $\langle\kappa,g\rangle$ has infinite abelianization, $\epsilon=1$.
After replacing $g$ by $gx$, if necessary,
we may assume that $\eta=1$.
Hence $gx^2=x^2g$ and $gy=yg$.
If $o(g)=2$ then $g^2$ is in the centre of $\kappa$, and so $g^2=x^{2m}$, for some $m\not=0$.
The exponent $s$ must then be odd, for otherwise $\pi$ would have non-trivial 2-torsion:
$(gy^{\frac{s}2}x^{-\frac{m}2})^2=1$.

If $\beta$ is generated by elements of finite order then the above argument 
implies that $\beta$ acts on $\kappa$ through a fixed subgroup of order 2 in
$Out(\kappa)\cong(\mathbb{Z}/2\mathbb{Z})^2$,
and that $\pi_1(F^+)=\kappa^+=\langle{x^2,y}\rangle$ is central in $\pi^+=\pi_1(M^+)$.
Therefore either $B$ is orientable or $B=Kb$, 
$\mathbb{A}$, $\mathbb{M}b$ or $P(2,2)$.
Moreover, $\beta$ has no element of order 4.
For if $o(g)=4$, $gxg^{-1}=xy^s$ and $gy=yg$ then $g^2xg^{-2}=xy^{2s}$.
If we set $h=g^2y^{-s}$ then $o(h)=2$, $hx=xh$ and $hy=yh$,
contrary to the earlier argument. 
Hence $B$ cannot be $S(2,4,4)$.

If $B=\mathbb{A}$ or $\mathbb{M}b$ then $\pi^{orb}(B)$ is a semidirect product
$D_\infty\rtimes\mathbb{Z}$, and so $\pi\cong{G}\rtimes\mathbb{Z}$,
where $G$ is an extension of $D_\infty$ by $\kappa$.
Thus $G\cong{J*_\kappa{K}}$, where  $[J:\kappa]=[K:\kappa]=2$.
Since $G$ is torsion-free $J\cong{K}\cong\kappa$.
We then find that $G$ is virtually $\mathbb{Z}^3$ and
$G/G'\cong\mathbb{Z}\oplus(\mathbb{Z}/2\mathbb{Z})^2$.
Hence $G$ is the non-orientable flat 3-manifold group $B_3$.
But $Out(B_3)$ is finite (see \S3 of Chapter 8 of \cite{Hi02}) and so $\pi$
must be virtually $\mathbb{Z}^4$.
Thus there are no such examples. 
\end{proof}

We shall show next that all the bases allowed by this theorem are realized by such Seifert fibred 4-manifolds.

Products of $\mathbb{N}il^3$-manifolds with $S^1$ give examples realizing the
five possibilities with $F=T$ and $\pi_1(F)$ central.
In these cases the fibration is the canonical one,
and is the orbit space projection of a $T$-action.

If $N$ is a $\mathbb{N}il^3$-manifold with Seifert base $B$ and
such that the general fibre has non-trivial image in $H_1(N;\mathbb{F}_2)$
then we may construct examples $M$ with base $B$
and general fibre $F=Kb$ as total spaces of $S^1$-bundles over $N$.
There are such 3-manifolds $N$ with base $B=T$, $Kb$ or $S(3,3,3)$.
The simplest are $Kb$-bundles over $T$, 
with groups
\[
\langle{a,b,x,y}\mid [a,b]=x^{2q}\!,~{a,b}\leftrightharpoons{x,y},~xyx^{-1}=y^{-1}\rangle,
\]
and $Kb$-bundles over $Kb$, with groups
\[
\langle{a,b,x,y}\mid{aba^{-1}b=x^{2q}\!,}~ax=xa,~bxb^{-1}=x^{-1}\!,~{a,b}\leftrightharpoons{y},~
xyx^{-1}=y^{-1}\rangle,
\]
where $q\not=0$. The groups
\[
\langle{c,u,x,y}\mid{c^3=u^3=x^2\!,}~(cu)^3=x^{4-2q}\!,
~{c,u}\leftrightharpoons{x,y},~xyx^{-1}=y^{-1}\rangle,
\]
are the groups of $\mathbb{N}il^3\times\mathbb{E}^1$-manifolds
with base $S(3,3,3)$ and general fibre $Kb$ if $q\not=0$.

There are also examples with base $B=S(2,2,2,2)$ or $S(2,3,6)$.
For instance, the groups
\[
\langle{j,a,b,x,y}\mid{j^2=y,} ~jaj^{-1}=a^{-1}\!, ~jbj^{-1}=b^{-1}\!,~ jxj^{-1}=xy^{-1}\!, 
\]
\[
[a,b]=x^{2q}\!,~a,b\leftrightharpoons{x,y},~xyx^{-1}=y^{-1}\rangle
\]
and
 \[
\langle{k,a,b,x,y}\mid{k^6=x^2y^3\!,}~kak^{-1}=b,~kbk^{-1}=a^{-1}b, ~kxk^{-1}=xy, ky=yk,
\]
\[
[a,b]=x^{2q}\!,~a,b\leftrightharpoons{x,y},~xyx^{-1}=y^{-1}\rangle
\]
are the groups of $\mathbb{N}il^3\times\mathbb{E}^1$-manifolds
with base $S(2,2,2,2)$ or $S(2,3,6)$ (respectively) and general fibre $Kb$ if $q\not=0$.
Finally, the groups
\[
\langle{d,j,x,y}\mid(jd^2)^2=x^{2q}y^{-1\!},~j^2=y^{-1}\!,~dx=xd,~dy=yd,~jxj^{-1}=xy, 
\]
\[
xyx^{-1}=y^{-1}\rangle
\]
are the groups of $\mathbb{N}il^3\times\mathbb{E}^1$-manifolds
with base $P(2,2)$ and general fibre $Kb$ if $q\not=0$.

In each case the constraint $q\not=0$ ensures that the group is not virtually $\mathbb{Z}^4$,
and the Seifert fibration is the canonical one.

For the remaining two geometries the general fibre is $T$ 
and the action of $\pi^{orb}(B)$ on the fundamental group of the fibre has 
infinite image in $Aut(\pi_1(T))\cong{GL(2,\mathbb{Z})}$,
since $M$ is neither flat nor a  $\mathbb{N}il^3\times\mathbb{E}^1$-manifold.
Since $\pi^{orb}(B)$ is solvable and $GL(2,\mathbb{Z})$ is virtually free this image
is virtually $\mathbb{Z}$, and so $B$  must fibre over $S^1$ or $\mathbb{I}$.
If the geometry is $\mathbb{S}ol^3\times\mathbb{E}^1$ then
the general fibre in the orbifold fibration of $B$ is $S^1$,
by Theorem 4 of \cite{Hi15}.
Hence $B$ is one  of $T$,  $Kb$, $\mathbb{A}$, $\mathbb{M}b$, 
$S(2,2,2,2)$, $P(2,2)$ or $\mathbb{D}(2,2)$.
Each of these possibilities occurs \cite{Hi15}.
Inspection of the tables on pages 219--230 of \cite{De} shows that the same seven bases 
are the bases realized by  $\mathbb{N}il^4$-manifolds.

\end{document}